%% file: weiss_equa.tex
\documentclass{amuc}
\usepackage{color,graphicx,wrapfig}
\begin{document}

\title[SHS - high activation energy]{Self-propagating High temperature Synthesis (SHS) in the
high activation energy regime}


\author{R. Monneau}
\address{Ecole Nationale des Ponts et Chauss\'ees, CERMICS,
 6 et 8 avenue Blaise Pascal,
Cit\'e Descartes Champs-sur-Marne, 77455 Marne-la-Vall\'ee  Cedex 2, France}
\email{monneau@cermics.enpc.fr}
\urladdr{http://cermics.enpc.fr/\textasciitilde monneau/home.html}


\author{G.S. Weiss}
\address{Graduate School of Mathematical Sciences,
University of Tokyo,
3-8-1 Komaba, Meguro, Tokyo,
153-8914 Japan}
\curraddr{until Sept. 2005 guest of the Max Planck Institute
for Mathematics in the Sciences,
Inselstr. 22, D-04103 Leipzig, German}
\email{gw@ms.u-tokyo.ac.jp}
\urladdr{http://www.ms.u-tokyo.ac.jp/\textasciitilde gw}


\thanks{G.S. Weiss has been partially supported by the Grant-in-Aid
15740100 of the Japanese Ministry of Education and partially supported
by a fellowship of the Max Planck Society. Both authors thank the Max Planck
Institute for Mathematics in the Sciences for the hospitality
during their stay in Leipzig.}




{}

\begin{abstract}
We derive the precise limit of SHS in the high
activation energy scaling suggested by B.J. Matkowksy-G.I. Sivashinsky 
in 1978 and by 
A. Bayliss-B.J. Matkowksy-A.P. Aldushin in 2002. 
In the time-increasing case the limit turns out to be
the Stefan problem for supercooled water with spatially
inhomogeneous coefficients.\\
Although the present paper leaves open mathematical
questions concerning the convergence, our precise form
of the limit problem suggest a strikingly simple explanation
for the numerically observed pulsating waves.  
\end{abstract}

\maketitle

\def\esssup{\hbox{\rm esssup }}
\def\R{\hbox{\bf R}}
\def\N{\hbox{\bf N}}
\def\sr{\hbox{\small\bf R}}
\def\sn{\hbox{\small\bf N}}
\def\supp{\hbox{\rm supp }}
\newtheorem{theo}{\textsc{Theorem}}[section]
\newtheorem{lem}[theo]{\textsc{Lemma}}
\newtheorem{pro}[theo]{\textsc{Proposition}}
\newtheorem{cor}[theo]{\textsc{Corollary}}
\newtheorem{defi}[theo]{\textsc{Definition}}
\newtheorem{rem}[theo]{\textsc{Remark}}
\newtheorem{hyp}[theo]{\textsc{Hypothesis}}
\newtheorem{pers}[theo]{\textsc{Perspectives}}
\newtheorem{conj}[theo]{\textsc{Conjecture}}

\section{Introduction}
The system
\begin{equation}
\begin{array}{l}\label{solid}
\partial_t u - \Delta u = v f(u)\\
\partial_t v = -v f(u)\; ,
\end{array}
\end{equation}
where $u$ is the normalized temperature, $v$ is the
normalized concentration of the reactant
and the non-negative nonlinearity $f$ describes the
reaction kinetics,
is a simple but widely used model for solid combustion
(i.e. the case of the Lewis number being $+\infty$).
In particular it is being used to model the industrial
process of Self-propagating High temperature Synthesis (SHS).
In the case of high activation energy interesting phenomena
like the instability of planar waves, fingering and 
helical waves are observed.\\
Since the seventies (and possibly even earlier)
it has been argued that the problem is for high
activation energy related to a Stefan problem describing
the freezing of supercooled water (see \cite{1978}, \cite[p. 57]{frankel}).
In \cite{1978} B.J. Matkowsky-G.I. Sivashinsky
derived a formal singular limit containing a jump condition
for the temperature on the interface.
Later the Stefan problem for supercooled water -- the
intuitive limit -- became the basis
for numerous papers focusing on stability analysis of
(\ref{solid}), fingering, helical waves etc.
(see for example \cite{frankel1},\cite{frankel2},\cite{frankel3},\cite{frankel4},\cite{frankel5},\cite{frankel6},\cite{frankel7},\cite{aldushin1},\cite{hotspots}).
\\
Surprisingly there are few {\em mathematical} results on
the subject:
In \cite{loubeau} E. Logak-V. Loubeau proved existence of a planar
wave in one-space dimension and gave a rigorous proof
for convergence as the activation energy goes to infinity.
\\
Instability of the planar wave
for a special linearization
(and high activation energy) is due to \cite{bonnet}.
\\[.5cm]
In the present paper we argue that the SHS system
converges to the irreversible
Stefan problem for supercooled water.
As the initial data of the reactant concentration
enters the equation as the activation energy goes to
infinity, our result also suggests a {\em surprisingly 
simple explanation
for the numerically observed pulsating
waves} (cf. \cite{aldushin1} and \cite{hotspots}),
namely that they are caused by the spatial inhomogeneity
$v^0$ (or $Y^0$, respectively) in the below equation
and are therefore mathematically related to
the pulsating waves in \cite{hamel}.\\
In the time-increasing case we give a rigorous
convergence proof in higher dimensions.
For general initial data in one space-dimension see
our forthcoming paper \cite{siam}.
\\[.5cm]
In the original setting by B.J. Matkowsky-G.I. Sivashinsky \cite[equation (2)]{1978},
\begin{equation}
\begin{array}{l}
\partial_t u_N - \Delta u_N = (1-\sigma_N)Ne^N v_N \exp(-N/u_N),\\
\partial_t v_N = -Ne^Nv_N \exp(-N/u_N),
\end{array}
\end{equation}
each limit $u_\infty$ of $u_N>0$ as $N\to\infty$
satisfies for $(\sigma_N)_{N\in \N} \subset  \subset [0,1)$
(for $\sigma_N \uparrow 1, N\to\infty$ the limit in this scaling
is the solution of the heat equation; 
cf. Section \ref{appl1}
and Theorem \ref{main})
\begin{equation}
\partial_t u_\infty-v^0 \partial_t \chi=\Delta u_\infty \hbox{ in } (0,+\infty)\times\Omega,
\end{equation}
where $v^0$ are the initial data of $v_\infty$ and
\[ \chi(t,x)\left\{\begin{array}{ll}
\in [0,1],& \esssup_{(0,t)} u_\infty(\cdot,x)\le 1\; ,\\
=1,& \esssup_{(0,t)} u_\infty(\cdot,x)>1\; ,\end{array}\right.\]
and in the time-increasing case,
\[ \chi(t,x)\left\{\begin{array}{ll}
=0,& u_\infty(t,x)<1\; ,\\
\in [0,1],& u_\infty(t,x)=1\; ,\\
=1,& u_\infty(t,x)>1\; .\end{array}\right.\]

In the SHS system with another scaling and a temperature threshold (see \cite[p. 109-110]{hotspots}),
\begin{equation}
\begin{array}{l}
\partial_t \theta_N - \Delta \theta_N\\
= (1-\sigma_N)
NY_N \exp((N(1-\sigma_N)(\theta_N-1))/
(\sigma_N + (1-\sigma_N)\theta_N))\chi_{\{\theta_N>\bar \theta\}},\\
\partial_t Y_N = -(1-\sigma_N)
NY_N \exp((N(1-\sigma_N)(\theta_N-1))/
(\sigma_N + (1-\sigma_N)\theta_N))\chi_{\{\theta_N>\bar \theta\}}
\end{array}
\end{equation}
where 
$N(1-\sigma_N)>>1, \sigma_N\in (0,1)$ and $\bar\theta \in (0,1)$,
each limit $\theta_\infty$ of $\theta_N$ satisfies (cf. Section \ref{appl2} and Theorem \ref{main})
\begin{equation}
\partial_t \theta_\infty-Y^0 \partial_t \chi=\Delta \theta_\infty \hbox{ in } (0,+\infty)\times\Omega,
\end{equation}
where $Y^0$ are the initial data of $Y_\infty$ and
\[ \chi(t,x)\left\{\begin{array}{ll}
\in [0,1],& \esssup_{(0,t)} \theta_\infty(\cdot,x)\le 1\; ,\\
=1,& \esssup_{(0,t)} \theta_\infty(\cdot,x)>1\; ,\end{array}\right.\]
and in the time-increasing case,
\[ \chi(t,x)\left\{\begin{array}{ll}
=0,& \theta_\infty(t,x)<1\; ,\\
\in [0,1],& \theta_\infty(t,x)=1\; ,\\
=1,& \theta_\infty(t,x)>1\; .\end{array}\right.\]
To our knowledge this precise form of the limit problem, i.e. 
the equation with the discontinuous
hysteresis term,
has not been known. Even in the time-increasing case it does
not coincide with the formal result in \cite{1978}.
\\[.5cm]
In the case that $\theta_\infty$ (or $u_\infty$, respectively) is increasing
in time
and $v^0$ (or $Y^0$, respectively) is constant, 
our limit problem coincides with the
Stefan problem for supercooled water,
an extensively studied ill-posed problem (for a survey see \cite{dewynne}). 
As it is a forward-backward parabolic equation it is not
clear whether one should expect uniqueness (see \cite[Remark 7.2]{friedman} for an example
of non-uniqueness in a related problem).\\
On the positive side,
much more is known about the Stefan problem for supercooled water
than the SHS system, e.g.
existence of a finger (\cite{ivantsov}), instability of the finger (\cite{langer}),
one-phase solutions (\cite{friedman});
those results, when combined with our convergence result,
suggest that similar properties should be true for the SHS
system.\\
It is interesting to observe that even in the 
time-increasing case our singular limit
{\em selects certain solutions} of the
Stefan problem for supercooled water.
For example, $u(t) =
(\kappa-1) \chi_{\{ t<1\}}+ \kappa \chi_{\{ t> 1\}}$
is for each $\kappa \in (0,1)$ a perfectly valid
solution of the Stefan problem for supercooled water,
but, as easily verified, it cannot be obtained from the ODE
\[ \partial_t u_\epsilon(t) = -\partial_t\exp(-{1\over \epsilon}\int_0^t \exp((1-1/(u_\epsilon(s)+1))/\epsilon)\> ds)\; \hbox{ as } \epsilon\to 0\; .\]
\section{Notation}
Throughout this article $\R^n$ will be equipped with the Euclidean
inner product $x\cdot y$ and the induced norm $\vert x \vert\> .$
$B_r(x)$ will denote the open $n$-dimensional ball of center
$x\> ,$ radius $r$ and volume $r^n\> \omega_n\> .$
When the center is not specified, it is assumed to be $0.$\\
When considering a set $A\> ,$ $\chi_A$ shall stand for
the characteristic function of $A\> ,$
while
$\nu$ shall typically denote the outward
normal to a given boundary.
The operator $\partial_t$ will mean the partial derivative
of a function in the time direction, $\Delta$ the
Laplacian in the space variables and ${\mathcal L}^n$ the $n$-dimensional Lebesgue measure.\\
Finally ${\bf W}^{2,1}_p$ denotes the parabolic Sobolev space
as defined in \cite{nina}.
\section{Preliminaries}
In what follows, $\Omega$ is a bounded $C^1$-domain in $\R^n$ and
\[ u_\epsilon\in \bigcap_{T \in (0,+\infty)}
{\bf W}^{2,1}_2((0,T)\times \Omega)\] is a strong solution
of the equation
\begin{equation}\label{sys}
\begin{array}{c}
\partial_t u_\epsilon(t,x) - \Delta u_\epsilon(t,x) = -v^0_\epsilon(x)\partial_t\exp(-{1\over \epsilon}\int_0^t g_\epsilon(u_\epsilon(s,x))\> ds)\; ,\\
u_\epsilon(0,\cdot)=u_\epsilon^0 \hbox{ in } \Omega, \nabla u_\epsilon\cdot \nu=0 \hbox{ on }
(0,+\infty)\times \partial\Omega \; ;
\end{array}
\end{equation}
here $g_\epsilon$ is a non-negative function
on $\R$ satisfying:\\
0) $g_\epsilon$ is for each $\epsilon\in (0,1)$
piecewise continuous with only one possible
jump at $z_0$, $g_\epsilon(z_0-)=g_\epsilon(z_0)=0$ in case of a jump, and
$g_\epsilon$ satisfies for each $\epsilon\in (0,1)$ and for every $z\in \R$ 
the bound $g_\epsilon(z)\le C_\epsilon (1+\vert z \vert)$.\\
1) $g_\epsilon/\epsilon \to 0$ as $\epsilon\to 0$
on each compact subset of $(-\infty,0)$.\\
2) for each compact subset $K$ of $(0,+\infty)$
there is $c_K>0$ such that 
$\min(g_\epsilon,c_K) \to c_K$ uniformly on $K$ as $\epsilon\to 0$.\\
The initial data satisfy $0\le v^0_\epsilon\le C<+\infty$, $v^0_\epsilon$ converges in $L^1(\Omega)$ to
$v^0$ as $\epsilon\to 0$,
$(u_\epsilon^0)_{\epsilon\in (0,1)}$ is bounded in $L^{2}(\Omega)$,
it is uniformly bounded from below by a constant $u_{\rm min}$,
and it converges in $L^1(\Omega)$ to $u^0$ as $\epsilon\to 0$.
\begin{rem}
Assumption 0) guarantees existence of a global strong solution for each $\epsilon\in (0,1)$.
\end{rem}
\section{The High Activation Energy Limit}
\begin{theo}\label{main}
The family $(u_\epsilon)_{\epsilon\in (0,1)}$ is for each
$T\in (0,+\infty)$ precompact in $L^1((0,T)\times \Omega)$, and
each limit $u$ of $(u_\epsilon)_{\epsilon\in (0,1)}$ as a sequence
$\epsilon_m\to 0,$ satisfies in the sense of distributions the 
initial-boundary value problem
\begin{equation}\label{stefan}
\partial_t u-v^0 \partial_t \chi =\Delta u \hbox{ in } (0,+\infty)\times\Omega,
\end{equation}
\[ u(0,\cdot)=u^0+v^0 H(u^0) \hbox{ in } \Omega\; , \; \nabla u\cdot \nu=0 \hbox{ on }
(0,+\infty)\times \partial\Omega\; ,\]
\[ \hbox{where } \chi(t,x) \left\{\begin{array}{ll}
\in [0,1],& \esssup_{(0,t)} u(\cdot,x)\le 0\; ,\\
=1,& \esssup_{(0,t)} u(\cdot,x)>0\; ,\end{array}\right.\]
and $H$ is the maximal monotone graph
\[ H(z) \left\{\begin{array}{ll}
=0,& z<0,\\ 
\in [0,1],& z=0,\\
=1,& z>0\; .\end{array}\right.\]
Moreover, $\chi$ is increasing in time and
$u$ is a supercaloric function.\\
If $(u_\epsilon)_{\epsilon\in (0,1)}$ satisfies
$\partial_t u_\epsilon\ge 0$ in $(0,T)\times \Omega$,
then $u$ is a solution of the Stefan problem for
supercooled water, i.e.
\[ \partial_t u-v^0 \partial_t H(u) =\Delta u \hbox{ in } (0,+\infty)\times\Omega\; .\]
\end{theo}
\begin{rem}
Note that assumption 1) is only needed to prove the second statement
``If ....''. 
\end{rem}
\proof
{\sl Step 0 (Uniform Bound from below):}\\
Since $u_\epsilon$ is supercaloric, it is bounded from
below by the constant $u_{\rm min}$.\\
{\sl Step 1 ($L^2((0,T)\times \Omega)$-Bound):}\\
The time-integrated function $v_\epsilon(t,x) := \int_0^t u_\epsilon(s,x)\> ds,$
satisfies
\begin{equation}
\partial_t v_\epsilon(t,x) - \Delta v_\epsilon(t,x) = w_\epsilon(t,x) + u_\epsilon^0(x)
\end{equation} where $w_\epsilon$ is a measurable function satisfying
$0\le w_\epsilon \le C.$
Consequently
\[ \int_0^T \int_\Omega (\partial_t v_\epsilon)^2 \; + \;
{1\over 2} \int_\Omega |\nabla v_\epsilon|^2(T)
\; = \; \int_0^T \int_\Omega (w_\epsilon + u_\epsilon^0) \partial_t v_\epsilon
\]\[ \le \; {1\over 2} \int_0^T \int_\Omega (\partial_t v_\epsilon)^2 \; + \;
{T\over 2} \int_\Omega (C+ |u_\epsilon^0|)^2\; ,\]
implying
\begin{equation}\label{est1}
\int_0^T \int_\Omega u_\epsilon^2 \; \le \;
T \int_\Omega (C+ |u_\epsilon^0|)^2\; . 
\end{equation}
{\sl Step 2 ($L^2((0,T)\times \Omega)$-Bound for $\nabla\min(u_\epsilon,M)$:}\\
For 
\[
G_M(z) := \left\{ \begin{array}{l}
z^2/2, z<M,\\
Mz-M^2/2, z\ge M\; ,\end{array}\right.\]
and
any $M\in\N,$
\[ \int_\Omega G_M(u_\epsilon)-G_M(u^0_\epsilon)
\; + \; \int_0^T \int_\Omega |\nabla \min(u_\epsilon,M)|^2
\]\[ = \; \int_0^T \int_\Omega -v^0_\epsilon \min(u_\epsilon,M)
\partial_t \exp(-{1\over \epsilon}\int_0^t g_\epsilon(u_\epsilon(s,x))\> ds)\; .\]
As $\partial_t \exp(-{1\over \epsilon}\int_0^t g_\epsilon(u_\epsilon(s,x))\> ds)\le 0,$
we know that $\partial_t \exp(-{1\over \epsilon}\int_0^t g_\epsilon(u_\epsilon(s,x))\> ds)$
is bounded in $L^\infty(\Omega;L^1((0,T))),$ and
\[ \int_0^T \int_\Omega -v^0_\epsilon \min(u_\epsilon,M)
\partial_t \exp(-{1\over \epsilon}\int_0^t g_\epsilon(u_\epsilon(s,x))\> ds)\]\[
\le \; C\> \int_\Omega \sup_{(0,T)} \max(\min(u_\epsilon,M),0) \;\le \;
CM {\mathcal L}^n(\Omega)\; .\]
\\
{\sl Step 3 (Compactness):}
\\
Let $\chi_M:\R\to\R$ be a smooth non-increasing function satisfying
$\chi_{(-\infty,M-1)}\le \chi_M \le \chi_{(-\infty,M)}$
and let $\Phi_M$ be the primitive such that $\Phi_M(z)=z$ for
$z\le M-1$ and $\Phi_M\le M$.
Moreover, let $(\phi_\delta)_{\delta\in (0,1)}$ be a family of mollifiers,
i.e.
$\phi_\delta\in C^{0,1}_0(\R^n;[0,+\infty))$ such that $\int \phi_\delta = 1$
and $\supp \phi_\delta \subset B_\delta(0)\> .$
Then, if we extend $u_\epsilon$ and $v^0_\epsilon$ by the value $0$ to the whole of $(0,+\infty)\times
\R^n,$ we obtain by the homogeneous Neumann data of $u_\epsilon$ that
\[ \partial_t \left(\Phi_M(u_\epsilon)*\phi_\delta\right)(t,x)
\]\[ = \; \left(\left(\chi_M(u_\epsilon)\left(\chi_\Omega\Delta u_\epsilon\> - \> v^0_\epsilon
\partial_t \exp(-{1\over \epsilon}\int_0^t g_\epsilon(u_\epsilon(s,x))\> ds)\right)
\right)*\phi_\delta\right)(t,x)
\]\[ = \; \int_{\sr^n}\chi_M(u_\epsilon)(t,y)
\bigg( \chi_\Omega(y)\Delta u_\epsilon(t,y)\]\[ - \> v^0_\epsilon(y)
\partial_t \exp(-{1\over \epsilon}\int_0^t g_\epsilon(u_\epsilon(s,y))\> ds)\bigg)
\phi_\delta(x-y)\> dy
\]\[ =\;\int_{\sr^n}\phi_\delta(x-y)\bigg(
-\chi_M'(u_\epsilon(t,y))\chi_\Omega(y)|\nabla u_\epsilon(t,y)|^2 \]\[- \>
\chi_M(u_\epsilon(t,y)) v^0_\epsilon(y)
\partial_t \exp(-{1\over \epsilon}\int_0^t g_\epsilon(u_\epsilon(s,y))\> ds)\bigg)\]\[
+ \>\chi_M(u_\epsilon(t,y)) \chi_\Omega(y)\nabla u_\epsilon(t,y)
\cdot \nabla \phi_\delta(x-y)\> dy\; .\]
Consequently 
\[ \int_0^T \int_{\sr^n} |\partial_t \left(\Phi_M(u_\epsilon)*\phi_\delta\right)|
\; \le \; C_1(\Omega,C,M,\delta,T)\]
and 
\[ \int_0^T \int_{\sr^n} |\nabla \left(\Phi_M(u_\epsilon)*\phi_\delta\right)|
\; \le \; C_2(\Omega,M,\delta,T)\; .\]
It follows that $(\Phi_M(u_\epsilon)*\phi_\delta)_{\epsilon\in (0,1)}$
is for each $(M,\delta,T)$ precompact in $L^1((0,T)\times \R^n)$.
\\
On the other hand
\[ \int_0^T \int_{\sr^n} |\Phi_M(u_\epsilon)*\phi_\delta-\Phi_M(u_\epsilon)|
\; \le \; C_3 \left(\delta^2\int_0^T \int_\Omega |\nabla \Phi_M(u_\epsilon)|^2\right)^{1\over 2}
\]\[ +\; 2(M-u_{\rm min})T{\mathcal L}^n(B_{\delta}(\partial\Omega))
\;\le \; C_4 (C,\Omega,u_{\rm min},M,T)\> \delta\; .\]
Combining this estimate
with the precompactness of $(\Phi_M(u_\epsilon)*\phi_\delta)_{\epsilon\in (0,1)}$ we obtain that $\Phi_M(u_\epsilon)$ is for
each $(M,T)$ precompact in $L^1((0,T)\times \R^n)$.
Thus, by a diagonal sequence
argument, we may take a sequence $\epsilon_m\to 0$ such that
$\Phi_M(u_{\epsilon_m})\to z_M$ a.e. in $(0,+\infty)\times \R^n$
as $m\to\infty$, for every $M\in \N.$
At a.e. point of the set $\{ z_M < M-1\}, u_{\epsilon_m}$
converges to $z_M.$ At each point $(t,x)$ of the remainder $\bigcap_{M\in \sn} 
\{ z_M \ge M-1\},$ the value $u_{\epsilon_m}(t,x)$ must for large $m$
(depending on $(M,t,x)$) be larger than $M-2$. But that means that
on the set $\bigcap_{M\in \sn} 
\{ z_M \ge M-1\},$ the sequence $(u_{\epsilon_m})_{m\in \sn}$ converges 
a.e. to $+\infty.$ It follows that $(u_{\epsilon_m})_{m\in \sn}$
converges a.e. in $(0,+\infty)\times \Omega$ to a function 
$z:(0,+\infty)\times\Omega\to \R\cup\{ +\infty\}$.
But then, as $(u_{\epsilon_m})_{m\in \sn}$ is for each $T\in (0,+\infty)$
bounded in $L^2((0,T)\times\Omega)$, $(u_{\epsilon_m})_{m\in \sn}$
converges by Vitali's theorem (stating that a.e. convergence
and a non-concentration condition in $L^p$ imply 
in bounded domains $L^p$-convergence) for each $p\in [1,2)$ in $L^p((0,T)\times\Omega)$
to the weak $L^2$-limit $u$ of $(u_{\epsilon_m})_{m\in \sn}$.
It follows that 
${\mathcal L}^{n+1}(\bigcap_{M\in \sn} 
\{ z_M \ge M-1\})={\mathcal L}^{n+1}(\{ u=+\infty\})=0$.
\\
{\sl Step 4 (Identification of the Limit Equation in $\esssup_{(0,t)} u>0$):}\\
Let us consider 
$(t,x)\in (0,+\infty)\times \Omega$ such that
$u_{\epsilon_m}(s,x)\to u(s,x)$ for a.e. $s \in (0,t)$
and $u(\cdot,x) \in L^2((0,t))$.
In the case $\esssup_{(0,t)} u(\cdot,x)>0$,
we obtain by Egorov's theorem and assumption 2) that
$\exp(-{1\over {\epsilon_m}}\int_0^t g_{\epsilon_m}(u_{\epsilon_m}(s,x))\> ds)\to 0$
as $m\to \infty$.
\\
{\sl Step 5 (The case $\partial_t u_\epsilon \ge 0$):}\\
Let $(t,x)$ be such that $u_{\epsilon_m}(t,x)\to u(t,x)=\lambda<0$: Then
by assumption 1),
\[ \exp(-{1\over {\epsilon_m}}\int_0^t g_{\epsilon_m}(u_{\epsilon_m}(s,x))\> ds)
\ge \exp(-t{\max_{[u_{\rm min},\lambda/2]} g_{\epsilon_m}\over {\epsilon_m}})
\to 1\hbox{ as }m\to \infty\; .\]
\begin{rem}
1) For a more general result in one space-dimension see the
forthcoming paper \cite{siam}.
\\
2) We also obtain a rigorous convergence result
in the case of (higher dimensional) traveling waves 
with suitable conditions at infinity. In this case
our $L^2(W^{1,2})$-estimate (Step 2) implies
a no-concentration property of the time-derivative. 
\end{rem}
\section{Applications}\label{applications}
Although the limit equation is an ill-posed problem,
the convergence to the limit seems to be robust
with respect to perturbations of the $\epsilon$-system and the scaling:
here we mention two examples of different systems leading to
the same limit. Other examples can be found in mathematical
biology (see \cite{kawasaki} and \cite{maini}). 
\subsection{The Matkowsky-Sivashinsky scaling}\label{appl1}
We apply our result to the scaling in \cite[equation (2)]{1978}, i.e.
\begin{equation}
\begin{array}{l}
\partial_t u_N - \Delta u_N = (1-\sigma_N)Nv_N \exp(N(1-1/u_N)),\\
\partial_t v_N = -Nv_N \exp(N(1-1/u_N)),
\end{array}
\end{equation}
where the normalized temperature $u_N$ and the normalized
concentration $v_N$ are non-negative,
$(\sigma_N)_{N\in \N} \subset  \subset [0,1)$ (in the case
$\sigma_N \uparrow 1, N\to\infty$ the limit equation
in the scaling as it is would be the heat equation,
but we could still apply our result to $u_N/(1-\sigma_N)$)
and the activation energy
$N\to \infty$.
\\
Setting $u_{\rm min} := -1, \epsilon := 1/N, u_\epsilon := u_N-1$ and
\[ g_\epsilon(z) := \left\{\begin{array}{ll}
\exp((1-1/(z+1))/\epsilon), z>-1\\
0, z\le -1\end{array}\right .\]
and integrating the equation for $v_N$ in time, 
we see that the assumptions of Theorem \ref{main} are satisfied
and we obtain
that each limit $u_\infty,\sigma_\infty$ of $u_N,\sigma_N$ satisfies
\begin{equation}\label{matsiv0}
\partial_t u_\infty-(1-\sigma_\infty)
v^0 \partial_t \chi=\Delta u_\infty \hbox{ in } (0,+\infty)\times\Omega,
\end{equation}
where
\[\chi(t,x)\left\{\begin{array}{ll}
\in [0,1],& \esssup_{(0,t)} u_\infty(\cdot,x)\le 1\; ,\\
=1,& \esssup_{(0,t)} u_\infty(\cdot,x)>1\; ,\end{array}\right.\]
and in the time-increasing case,
\begin{equation}\label{matsiv}
\partial_t u_\infty-(1-\sigma_\infty)
v^0 \partial_t H(u_\infty)=\Delta u_\infty \hbox{ in } (0,+\infty)\times\Omega,
\end{equation}
\[ u_\infty(0,\cdot)=u^0+v^0H(u^0) \hbox{ in } \Omega\; , \; \nabla u_\infty\cdot \nu=0 \hbox{ on }
(0,+\infty)\times \partial\Omega\; ,\]
where $v^0$ are the initial data of $v_\infty$.
Moreover, $\chi$ is increasing in time and
$u_\infty$ is a supercaloric function.
\subsection{SHS in another scaling with temperature threshold}\label{appl2}
Here we consider (cf. \cite[p. 109-110]{hotspots}), 
i.e.
\begin{equation}
\begin{array}{l}
\partial_t \theta_N - \Delta \theta_N \\
= (1-\sigma_N)
NY_N \exp((N(1-\sigma_N)(\theta_N-1))/
(\sigma_N + (1-\sigma_N)\theta_N))\chi_{\{\theta_N>\bar \theta\}},\\
\partial_t Y_N = -(1-\sigma_N)
NY_N \exp((N(1-\sigma_N)(\theta_N-1))/
(\sigma_N + (1-\sigma_N)\theta_N))\chi_{\{\theta_N>\bar \theta\}}
\end{array}
\end{equation}
where 
$N(1-\sigma_N)>>1, \sigma_N\in (0,1)$ and
the constant $\bar\theta \in (0,1)$ is a threshold parameter
at which the reaction sets in.
\\
Setting $u_{\rm  min}=-1, \epsilon := 1/(N(1-\sigma_N)), 
\kappa(\epsilon):= 1-\sigma_N, u_\epsilon := \theta_N-1$,
\[ g_\epsilon(z) := \left\{\begin{array}{ll}
\exp((z/(\kappa(\epsilon)z+1))/\epsilon), z>\bar\theta-1\\
0, z\le \bar\theta-1\end{array}\right .\]
and integrating the equation for $Y_N$ in time,
we see that the assumptions of Theorem \ref{main} are satisfied
and we obtain that each limit $u_\infty$ of $u_N$ satisfies
\begin{equation}\label{threshold0}
\partial_t u_\infty-v^0 \partial_t \chi=\Delta u_\infty \hbox{ in } (0,+\infty)\times\Omega,
\end{equation}
\[\chi(t,x)\left\{\begin{array}{ll}
\in [0,1],& \esssup_{(0,t)} u_\infty(\cdot,x)\le 1\; ,\\
=1,& \esssup_{(0,t)} u_\infty(\cdot,x)>1\; ,\end{array}\right.\]
and in the time-increasing case,
\begin{equation}\label{threshold}
\partial_t u_\infty-v^0 \partial_t H(u_\infty)=\Delta u_\infty \hbox{ in } (0,+\infty)\times\Omega,
\end{equation}
\[ u_\infty(0,\cdot)=u^0+v^0H(u^0) \hbox{ in } \Omega\; , \; \nabla u_\infty\cdot \nu=0 \hbox{ on }
(0,+\infty)\times \partial\Omega\; ,\]
where $v^0$ are the initial data of $v_\infty$.
Moreover, $\chi$ is increasing in time and
$u_\infty$ is a supercaloric function.
\section{Open questions}
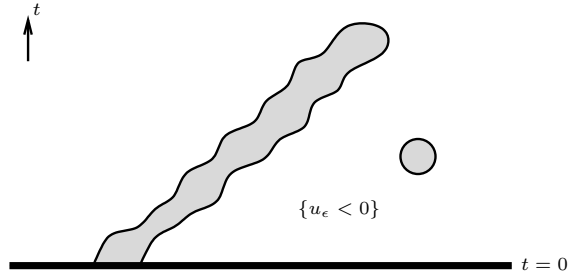
\begin{figure}
\begin{center}
\input{weiss_equa.pstex_t}
\end{center}
\caption{Is it possible for the solution to have a tiny peak
traveling at high speed?}\label{peak}
\end{figure}

The most pressing task is of course to 
study the existence or non-existence of ``peaking'' (cf. Figure
\ref{peak})
of the solution in the negative phase (for
the case of one space dimension see the forthcoming
paper \cite{siam}).
A related question
is whether $(u_\epsilon)_{\epsilon \in (0,1)}$ is bounded in $L^\infty$
in the case of uniformly bounded initial data. Although this
seems obvious, it is not obvious how to prevent concentration
close to the interface. \\
Uniqueness for the limit problem
(the irreversible Stefan problem for supercooled water) 
in general seems unlikely. One might however ask whether
time-global uniqueness holds in the case that $u$ is strictly increasing
in the $x_1$-direction. By the result in
\cite{hele} for the ill-posed Hele-Shaw problem, 
time-local uniqueness is likely to be true here, too.

\begin{acknowledgment}
We thank Stephan Luckhaus, Mayan Mimura, Stefan M\"uller, 
and Juan J.L. Vel\'azquez for discussions.
\end{acknowledgment}

\bibliographystyle{plain}
\bibliography{weiss_equa.bib}
\end{document}

%% file: weiss_equa.pstex_t
\begin{picture}(0,0)%
\includegraphics{weiss_equa.pstex}%
\end{picture}%
\setlength{\unitlength}{3947sp}%
\begingroup\makeatletter\ifx\SetFigFont\undefined%
\gdef\SetFigFont#1#2#3#4#5{%
  \reset@font\fontsize{#1}{#2pt}%
  \fontfamily{#3}\fontseries{#4}\fontshape{#5}%
  \selectfont}%
\fi\endgroup%
\begin{picture}(3508,1703)(136,-1030)
\put(1991,-654){\makebox(0,0)[lb]{\smash{{\SetFigFont{7}{8.4}{\rmdefault}{\mddefault}{\updefault}{\color[rgb]{0,0,0}$\{ u_\epsilon<0\}$}%
}}}}
\put(3393,-1004){\makebox(0,0)[lb]{\smash{{\SetFigFont{7}{8.4}{\rmdefault}{\mddefault}{\updefault}{\color[rgb]{0,0,0}$t=0$}%
}}}}
\put(326,602){\makebox(0,0)[lb]{\smash{{\SetFigFont{7}{8.4}{\rmdefault}{\mddefault}{\updefault}{\color[rgb]{0,0,0}$t$}%
}}}}
\end{picture}%